\documentclass[a4paper,12pt, reqno]{amsart}
\usepackage{hyperref}
\usepackage{a4wide}
\usepackage{amsthm}
\usepackage{amssymb}
\usepackage{longtable}

\newtheorem{thm}{Theorem}[section]

\newtheorem{lem}[thm]{Lemma}

\theoremstyle{definition}

\theoremstyle{remark}

\makeatletter
\let\c@equation\c@thm
\makeatother
\numberwithin{equation}{section}

\begin{document}
	%%% \topmatter
	
\title[A note on Euclidean cyclic cubic fields]{A note on Euclidean cyclic cubic fields}
\author{Srinivas Kotyada and Subramani Muthukrishnan}
\address{Institute of Mathematical Sciences, HBNI, CIT Campus, Taramani, Chennai 600 113, India}
\address{Chennai Mathematical Institute, SIPCOT IT Park, Siruseri, Chennai 603 103, India}
\email[Kotyada Srinivas]{srini@imsc.res.in}
\email[Subramani Muthukrishnan]{subramani@cmi.ac.in}

\begin{abstract}
Let $K$ be a cyclic cubic field and $\mathcal{O}_K$ be its ring of integers. In this note we prove that all cyclic cubic number fields with conductors in the interval $ [73,  11971]$ and with class number one are Euclidean. 
	
\end{abstract}
		
\maketitle
 
\section{Introduction}

We first recall that if $K$ is an algebraic number field and $\mathcal{O}_K$ its ring of integers, then  $\mathcal{O}_K$ is called Euclidean with respect to a given function $\phi : \mathcal{O}_K \to \mathbb{N} \cup \{0\}$ provided $\phi$ satisfies the following properties  
\begin{enumerate}
\item  $\phi (\alpha) = 0  \textrm{ if and only if } \alpha = 0$, and  
\item  for all $\alpha, \beta \not = 0 \in \mathcal{O}_K$ there exists a $\gamma \in \mathcal{O}_K$ such that $\phi(\alpha - \beta \gamma) < \phi(\beta).$ 
\end{enumerate}
Moreover, it is called norm-Euclidean if it is Euclidean with respect to the absolute value norm. 

\medskip

Now, let $K$ be a cyclic cubic field with discriminant $f^2,$ where $f$ is the conductor of $K.$ In 1969, J. R. Smith \cite{smith} proved that the cyclic cubic fields with conductors $7,9,...,67$ are norm-Euclidean. Further, in the same paper he showed that the fields with conductors $73, 79, 97, 139,151$ and $163<f<10^4$ are \emph{not} norm-Euclidean. The object of this note is to show that all cyclic cubic fields with conductors $f \in [73, 11971]$ are in fact Euclidean provided they have class number one.

\medskip

It is well known that if the conductor $f$ of the cyclic cubic field $K$ has $t$ distinct prime factors, then the class number of $K$ is divisible by $ 3^{t-1}$ (see appendix of Heilbronn's paper \cite{Heil} for a proof). Thus a cyclic cubic field with class number one must have prime conductor $f$. Moreover,  a necessary condition for a cyclic cubic field to have class number one is that its conductor is either $9$ or  a prime in the residue class $1 \pmod{6}$ ( see \cite{Heil}, \cite{G-S} ). Accordingly, from now onwards, we shall be dealing with only those cyclic cubic fields $K$ with prime conductor $f$ satisfying $f \equiv  1 \pmod{6}$. Our main aim in this note is to prove the following

\begin{thm}
Let $K$ be a cyclic cubic field with conductor $f,$ satisfying $ 73 \leq f \leq 11971$ and let  $\mathcal{O}_K$ be its ring of integers. Then $\mathcal{O}_K$ is Euclidean if and only if it has class number one. 
\end{thm}

\medskip

\noindent
\textbf{Proof.} It is easy to show that if a number ring $\mathcal{O}_K$ is Euclidean, then it has class number one.  To prove the converse,
we use a result of Harper and Ram Murty (Theorem \ref{thm:ram} below) which gives a useful criteria to establish the Euclidean algorithm for certain number fields. In order to state their theorem, we need to define the concept of a set of \emph{admissible} primes (see \cite{Clark} for details) in a number ring, which we do below.
	
\medskip

Assume that $\mathcal{O}_K$ has class number one.  Let $\pi_1, \dots, \pi_s \in \mathcal{O}_K$ be distinct non-associate primes.  A set of primes $\{\pi_1, \dots , \pi_s\}$ is called an \emph{admissible} set of primes if, for all $\beta = \pi_1^{a_1} \dots \pi_s^{a_s}$ with $a_i$ non-negative integers, every co-prime residue class$\pmod{\beta}$ can be represented by a unit $\varepsilon \in \mathcal{O}_K^\times.$ In other words, the  set $\{\pi_1, \dots , \pi_s\}$ is \emph{admissible} if the canonical map $\mathcal{O}_K^\times \to \big(\mathcal{O}_K/(\pi_1^{a_1} \dots \pi_s^{a_s})\big)^\times$ is surjective. In fact, in \cite{Clark}, the authors showed that it is enough to take $a_1=a_2= \dots = a_s =2$ in the above definition, (i.e,. the set $\{\pi_1, \dots , \pi_s\}$ is \emph{admissible} if the canonical map $\mathcal{O}_K^\times \to \big(\mathcal{O}_K/(\pi_1^2 \dots \pi_s^2)\big)^\times$ is a surjective).

\medskip

The precise statement of Harper and Ram Murty mentioned above is as follows:
\begin{thm} (M. Ram Murty, M. Harper \cite{Harper-Murty}) \label{thm:ram}
Let $K/\mathbb{Q}$ be abelian of degree $n$ with $\mathcal{O}_K$ having class number one, that contains a set of admissible primes with $s$ elements.  Let $r$ be the rank of the unit group.  If $r+s \geq 3,$ then $\mathcal{O}_K$ is Euclidean.
\end{thm}

The well known Dirichlet's unit theorem states that the rank $r$ of the group of units in $\mathcal{O}_K$ is given by $ r = r_1 + r_2 - 1$, where $r_1$ is the number of real embeddings and $r_2$ is the number of conjugate pairs of complex embeddings of $K$. In our case, $K$ is cyclic cubic field which means the Galois group over $\mathbb{Q}$ is cyclic of order three. This can only happen if $K$ is totally real. Thus, $ r = 3 - 1 = 2$. All we need to do now is to exhibit an \emph{admissible} set of primes with one element (i.e. $s =1$). For this we recall the following lemma from \cite{R-S-S}. Before stating the lemma, we recall the notion of non-Wieferich primes in a number field. 

\medskip

\noindent
 A prime $\mathfrak{p}$ in $\mathcal{O}_K$ is called Wieferich prime with respect to the base $\varepsilon \in \mathcal{O}_K^\times$ if
$$
\varepsilon^{N(\mathfrak{p})-1} \equiv 1 \pmod{\mathfrak{p}^2},
$$
where $N(.)$ is the absolute value norm.  If the above congruence does not hold for a prime $\mathfrak{p}$ in $\mathcal{O}_K$, then it is  called non-Wieferich prime to the base $\varepsilon.$

\begin{lem} \label{lem-first}
Let $\rho \in \mathcal{O}_K^\times$ be a unit and $\mathfrak{q}$ be an unramified prime ideal with odd prime norm $q$. If $\rho$ is a primitive root modulo $\mathfrak{q},$ and $\mathfrak{q}$ is a non-Wieferich prime to the base $\rho$, i.e., $\rho^{q-1} \not \equiv 1 \; (mod \; \mathfrak{q}^2)$, then $\rho$ generates the group $(\mathcal{O}_K/\mathfrak{q}^2)^\times.$
\end{lem}

\noindent
Thus, we need to find an unramified prime $\pi$ with odd prime norm such that the group $(\mathcal{O}_K/\pi)^\times$ has a primitive root $\varepsilon \in \mathcal{O}_K^\times$ and $\pi$ is a non-Wieferich prime with respect to $\varepsilon$. As the field $K$ is Galois of degree $3$ over $\mathbb{Q}$, this means that an unramified rational prime $p$ either splits completely  or remains as a prime in $\mathcal{O}_K$. It is well-known that a rational prime $p \textrm{ } (\not= f)$ splits completely in $K$ if and only if $p$ is a cube modulo $f$. By Euler's criterion, it follows that $p$ is a cube modulo $f$ if and only if $p^{\frac{f-1}{3}} \equiv 1 \pmod{f}$.

\medskip

In what follows, we shall exhibit a set of \emph{admissible} primes with one element for the cyclic cubic field $K$ with conductor $f = 73$. By Theorem \ref{thm:ram} and Lemma \ref{lem-first} it will then follow that this field is Euclidean as it has class number one. For all other fields in the range $ 73 < f \leq 11971$ with class number one, we shall give an algorithm which produces an \emph{admissible} set of primes with one element. The database \cite{J-R} gives all such fields. In the end, we shall list in a table the defining polynomial of all the class number one  cyclic cubic fields with conductors in the above range and the corresponding set of \emph{admissible} primes. This will complete the proof of our main theorem.

\noindent

Thus, we start with the cyclic cubic field $K$ with conductor $73$.  It is known that $K$ has class number one. 
%We note that $K = \mathbb{Q}(a)$ where $a$ is a root of the polynomial $x^3-x^2-24x+27$. 
The fundamental units $\varepsilon_1, \varepsilon_2$  for $K$  are 
\begin{align*} 
\varepsilon_1 := \frac{2}{3}a^2 - \frac{14}{3}a + 7, \\
\varepsilon_2 :=  \frac{4}{3}a^2 + \frac{14}{3}a - 7,
\end{align*}
where $a$ is a root of the defining polynomial $x^3-x^2-24x+27$ for $K$. This is obtained by using Sage programme.

\medskip

\noindent
Let us take the rational prime $p =3.$  It is unramified in $K$ since $p \nmid 73.$ 

\noindent
Also, as $3^{\frac{73-1}{3}} \equiv 1 \pmod {73}$ implies that $p$ splits completely in $\mathcal{O}_K$.
The prime ideal decomposition of $3$ in $\mathcal{O}_K$ is given by:
\begin{align*}
\left(3\right) = \left(\tfrac{1}{3}a^2 + \tfrac{2}{3}a - 11\right) \left(-\tfrac{1}{3}a^2 + \tfrac{1}{3}a + 6\right) \left(\tfrac{2}{3}a^2 + \tfrac{1}{3}a - 17\right).
\end{align*}

\noindent
Let $\pi$ be the prime element such that $(\pi) := \left(\frac{1}{3}a^2 + \frac{2}{3}a - 11\right).$ A simple calculation gives 
\begin{equation} \label{eqn-1}
\varepsilon_1 \equiv -1 \pmod{\pi},
\end{equation}
\noindent
and 
\begin{equation}\label{eqn-2}
\varepsilon_1^2 \equiv 16 \not\equiv 1 \pmod{{\pi}^2}.
\end{equation}
As $(\mathcal{O}_K/\pi)^\times$ has order 2, it follows from \eqref{eqn-1} that $\varepsilon_1$ is a primitive root modulo $\pi$. The equation \eqref{eqn-2} says that $\pi$ is a non-Wieferich prime with respect to $\varepsilon$. Thus, by Lemma \ref{lem-first} the set \{$\pi$\} is \emph{admissible}. Thus $K$ is Euclidean.

\medskip

\medskip

In section \ref{algo}, we present an algorithm to determine a set of \emph{admissible} primes with one element, section \ref{table} will contain table of admissible primes. We remark that, to save space, we only list a few fields which are Euclidean, the complete list (fields with conductors upto $11971$) is posted at
\begin{center}
\url{https://www.imsc.res.in/~srini/Euclidean-Cyclic-Cubic-Fields.txt} 
\end{center}

\medskip

\pagebreak

\section{Algorithm to find an admissible prime}\label{algo}

\medskip

\noindent
1. \hspace{.5cm} data $\leftarrow$ list of [defining polynomials, conductor] \\	
2. \hspace{.5cm} flat $\leftarrow$ false \\
3. \hspace{.5cm} foo $\leftarrow$ false\\ 
4. \hspace{.5cm} print ("Conductor|Defining Polynomial|Admissible prime") \\
5. \hspace{.5cm} for $u$ in data:\\
6. \hspace{1cm} prime $\leftarrow$	int(math.sqrt(u[1])) \\
7. \hspace{1cm} $K \leftarrow$ Number field with defining polynomial u[0] \\
8. \hspace{1cm} if (class number of $K$ is not $1$):\\
9. \hspace{1.5cm} break \\
10. \hspace{1cm} u $\leftarrow$ unit group associated with $K$\\
11. \hspace{1cm} eps $\leftarrow$ one of the fundamental units of u \\
12. \hspace{1cm} list p $\leftarrow$ [all cubic residues mod p] \\
13. \hspace{1cm} for q set of primes\\
14. \hspace{1.5cm} if q==2:\\
15. \hspace{2cm} continue\\
16. \hspace{1.5cm} if q$>$100000: \\
17. \hspace{2cm} break \\
18. \hspace{1.5cm} ((q modulo prime) in list p):\\
19. \hspace{2cm} f $\leftarrow$ prime decomposition of q in $K$\\
20. \hspace{2cm} for term in f:\\
21. \hspace{2.5cm} frc $\leftarrow$ fractional ideal in term \\
22. \hspace{2.5cm} for $d \in \{1,\dots q-1\}$\\
23. \hspace{3cm} if {(q-1 mod d is 0 and ($\textrm{eps}^d$) modulo frc ==1):} \\
24. \hspace{3.5cm} flat $\leftarrow$ true\\
25. \hspace{3.5cm} break\\
26. \hspace{2cm} if (flat is true):\\
27. \hspace{2.5cm} flat $\leftarrow$ false \\
28. \hspace{2.5cm} continue \\
29. \hspace{2cm} frcs $\leftarrow$ $(\textrm{frc})^2$\\
30. \hspace{2cm} if ($\textrm{eps}^{q-1}$ mod frcs is not 1) \\
31. \hspace{2.5cm} print(prime, polynomial, fractional ideal)\\
32. \hspace{3cm} foo $\leftarrow$ true \\
33. \hspace{3cm} break\\
34. \hspace{.5cm} foo is true:\\
35. \hspace{1cm} foo $\leftarrow$ false\\
36. \hspace{1cm} break\\
37. \hspace{.5cm} if(prime $> 100000$):\\
38. \hspace{1cm} break.

\pagebreak

\section{Table of admissible primes}\label{table}

\begin{center}
\begin{longtable}{|l|l|l|}
\caption{Cyclic Cubic number fields.} \label{tab:long} \\
         							
\hline \multicolumn{1}{|c|}{\textbf{$f$}} & \multicolumn{1}{c|}{\textbf{$polynomial$}} &
\multicolumn{1}{c|}{\textbf{$admissible$ $prime$}} \\ \hline 
\endfirsthead
		
\multicolumn{3}{c}%
{{\bfseries \tablename\ \thetable{} }} \\
\hline \multicolumn{1}{|c|}{\textbf{$f$}} & \multicolumn{1}{c|}{\textbf{$polynomial$}}&
\multicolumn{1}{c|}{\textbf{$admissible$ $prime$}} \\ \hline 
\endhead
		
\hline \multicolumn{3}{|r|}{{Continued on next page}} \\ \hline
\endfoot
		
\hline \hline
\endlastfoot
		
$79$ & $x^3-x^2-26x-41$ & $a+1$ \\
$97$ & $x^3-x^2-32x+79$ & $a^2+3a-21$ \\
$103$ & $x^3-x^2-34x+61$ & $\frac{2}{3}a^2-3a-\frac{8}{3}$ \\
$109$ & $x^3-x^2-36x+4$ &$\frac{35}{2}a^2+\frac{191}{2}a-12$ \\
$127$ & $x^3-x^2-42x-80$ & $9a^2+5a-153$ \\
$139$ & $x^3-x^2-46x-103$ & $a+4$ \\
$151$ & $x^3-x^2-50x+123$ & $a^2+3a-53$ \\
$157$ & $x^3-x^2-52x-64$  & $\frac{135}{2}a^2+\frac{977}{2}a+523$ \\
$181$ & $x^3-x^2-60x+67$ & $\frac{2}{5}a^2+\frac{21}{5}a-\frac{26}{5}$ \\
$193$ & $x^3-x^2-64x-143$ & $\frac{29}{3}a^2-65a-\frac{734}{3}$ \\
$199$ & $x^3-x^2-66x-59$ & $2a^2+16a+13$ \\
$211$ & $x^3-x^2-70x+125$ & $a^2-62$ \\
$223$ & $x^3-x^2-74x+256$ & $5313a^2+28133a-216059$ \\
$229$ & $x^3-x^2-76x+212$ & $266a^2+1704a-7627$ \\
$271$ & $x^3-x^2-90x-261$ & $3a^2-3a-152$ \\
$283$ & $x^3-x^2-94x-304$ & $448a^2+4698a+11855$ \\
$331$ & $x^3-x^2-110x+49$ & $\frac{3}{7}a^2+\frac{30}{7}a-1$ \\
$337$ & $x^3-x^2-112x-25$ & $4a^2-3a-428$ \\
$367$ & $x^3-x^2-112x-435$ & $\frac{1633}{3}a^2-\frac{13919}{3}a-3149$ \\
$373$ & $x^3-x^2-124x+221$ & $\frac{38}{7}a^2-\frac{475}{7}+\frac{733}{7}$ \\
$379$ & $x^3-x^2-126x-365$ & $\frac{13}{5}a^2+\frac{156}{5}a+74$ \\
$409$ & $x^3-x^2-136x+515$ & $\frac{19}{5}a^2-\frac{252}{5}a+134$ \\
$421$ & $x^3-x^2-140x+343$& $\frac{13}{7}a^2+\frac{134}{7}a-60$ \\
$433$ & $x^3-x^2-144x+16$ & $\frac{4257}{4}a^2-\frac{53501}{4}a+1472$ \\
$439$ & $x^3 - x^2 - 146x + 504$ & $(3, \frac{1}{6}a^2 + \frac{1}{6}a-16)$\\
$457$ & $x^3 - x^2 - 152x + 220$ & $(5, a - 2)$\\
$463$ & $x^3 - x^2 - 154x - 343$ & $(7, \frac{2}{7}a^2 - 23/7a - 30)$\\
$487 $ & $ x^3 - x^2 - 162x + 505$ & $ (5, a)$\\
$499 $ & $ x^3 - x^2 - 166x - 536 $ & $ (13, a + 6)$\\
$523 $ & $ x^3 - x^2 - 174x + 891 $ & $ (11, a)$\\
$541 $ & $ x^3 - x^2 - 180x - 521 $ & $ (7, -\frac{2}{7}a^2 - \frac{3}{7}a + \frac{251}{7})$\\
$571 $ & $ x^3 - x^2 - 190x + 719 $ & $ (7, -\frac{3}{7}a^2 + \frac{4}{7}a + \frac{354}{7})$\\
$577 $ & $ x^3 - x^2 - 192x - 171 $ & $ (3, \frac{1}{9}a^2 + \frac{11}{9}a - \frac{44}{3})$\\
$601 $ & $ x^3 - x^2 - 200x - 512 $ & $ (13, a + 3)$\\
$613 $ & $ x^3 - x^2 - 204x - 999 $ & $ (3, \frac{1}{3}a^2 - \frac{13}{3}a - 44)$\\
$619 $ & $ x^3 - x^2 - 206x - 321 $ & $ (3, a)$\\
$631 $ & $ x^3 - x^2 - 210x + 1075 $ & $ (43, a + 14)$\\
$643 $ & $ x^3 - x^2 - 214x + 1024 $ & $ (3, \frac{1}{6}a^2 +\frac{1}{2}a - \frac{68}{3})$\\
$661 $ & $ x^3 - x^2 - 220x + 1273 $ & $ (3, a + 1)$\\
$673 $ & $ x^3 - x^2 - 224x + 997 $ & $ (23, a + 3)$\\
$691 $ & $ x^3 - x^2 - 230x - 128 $ & $ (5, -\frac{1}{10}a^2 - \frac{23}{10}a + 74$\\
$727 $ & $ x^3 - x^2 - 242x - 1104 $ & $ (3, \frac{1}{6}a^2 + \frac{1}{6}a - 27)$\\
$733 $ & $ x^3 - x^2 - 244x - 1276 $ & $ (5, a - 1)$\\
$739 $ & $ x^3 - x^2 - 246x + 520 $ & $ (5, a)$\\
$751 $ & $ x^3 - x^2 - 250x - 1057 $ &  $ (7, -\frac{1}{7}a^2 + \frac{11}{7}a + 23)$\\
$757 $ & $x^3 - x^2 - 252x - 729 $ & $ (3, -\frac{1}{9}a^2 + \frac{10}{9}a + 19)$\\
$769 $ & $ x^3 - x^2 - 256x + 1481 $ & $ (5, \frac{1}{5}a^2 + a - \frac{166}{5})$\\
$787 $ & $ x^3 - x^2 - 262x + 991 $ & $ (31, a + 8)$\\
$811 $ & $ x^3 - x^2 - 270x - 1592 $ & $ (7, a - 2)$\\
$823 $ & $ x^3 - x^2 - 274x - 61 $ & $ (5, a + 2)$\\
$829 $ & $ x^3 - x^2 - 276x + 307 $ & $ (7, a + 2)$\\
$859 $ & $ x^3 - x^2 - 286x + 509 $ & $ (59, a + 28)$\\
$883 $ & $ x^3 - x^2 - 294x - 1439 $ & $ (17, a - 4)$\\
$907 $ & $ x^3 - x^2 - 302x + 739 $ & $ (11, -\frac{5}{11}a^2 - \frac{40}{11}a + \frac{1007}{11})$\\
$919 $ & $ x^3 - x^2 - 306x + 1872 $ & $ (29, a - 4)$\\
$967 $ & $ x^3 - x^2 - 322x - 1361 $ & $ (3, \frac{1}{9}a^2 - \frac{1}{3}a - \frac{226}{9})$\\
$991 $ & $ x^3 - x^2 - 330x + 2349 $ & $ (3, \frac{1}{3}a^2 + \frac{11}{3}a - 74)$\\
$997 $ & $ x^3 - x^2 - 332x + 480 $ & $ (5, a + 1)$\\
$1021 $ & $ x^3 - x^2 - 340x - 416 $ & $ (7, a - 2)$\\
$1033 $ & $ x^3 - x^2 - 344x - 1913 $ & $ (37, a - 1)$\\
$1039 $ & $ x^3 - x^2 - 346x - 2155 $ & $ (5, a)$\\
$1051 $ & $ x^3 - x^2 - 350x + 2608 $ & $ (53, a + 24)$\\
$1069 $ & $ x^3 - x^2 - 356x - 2336 $ & $ (29, a - 6)$\\
$1087 $ & $ x^3 - x^2 - 362x + 2335 $ & $ (5, a - 2)$\\
$1093 $ & $ x^3 - x^2 - 364x + 1012 $ & $ (3, \frac{1}{12}a^2 - \frac{1}{4}a - \frac{113}{6})$\\
$1117 $ & $ x^3 - x^2 - 372x - 2565 $ & $ (3, \frac{1}{3}a^2 - \frac{13}{3}a - 82)$\\
$1123 $ & $ x^3 - x^2 - 374x - 1331 $ & $ (5, a + 2)$\\
$1153 $ & $ x^3 - x^2 - 384x + 427 $ & $ (7, a)$\\
$1171 $ & $ x^3 - x^2 - 390x - 347 $ & $ (13, \frac{1}{13}a^2 + \frac{4}{13}a - \frac{214}{13})$\\
$1201 $ & $ x^3 - x^2 - 400x - 2491 $ & $ (13, a - 4)$\\
$1213 $ & $ x^3 - x^2 - 404x - 629 $ & $ (13, a - 6)$\\
$1231 $ & $ x^3 - x^2 - 410x + 1003 $ & $ (11, a + 4)$\\
$1237 $ & $ x^3 - x^2 - 412x - 1741 $ & $ (11, \frac{1}{11}a^2 - \frac{2}{11}a - \frac{245}{11})$\\
$1249 $ & $ x^3 - x^2 - 416x - 2313 $ & $ (3, a)$\\
$1279 $ & $ x^3 - x^2 - 426x + 2179 $ & $ (11, a + 3)$\\
$1291 $ & $ x^3 - x^2 - 430x + 3347 $ & $ (5, \frac{1}{5}a^2 + \frac{8}{5}a - \frac{278}{5})$\\

\end{longtable}
\end{center}

%\pagebreak

\section{Concluding remarks}

As and when new fields with class number one are added to the database \cite{J-R}, it is possible to check with the algorithm given in this paper to determine whether it is Euclidean or not. The Sage code for the algorithm written in this paper is posted at:
\begin{center}
	\url{https://www.imsc.res.in/~srini/Sage-Code-Cyclic-Cubic-Fields.txt}
\end{center}

\medskip

\noindent
\textbf{Acknowledgement:}
We sincerely thank Prof. M. Ram Murty for suggesting this problem and for some useful discussions. The second author would like to express his sincere gratitude to Prof. T.R. Ramadas for encouragement and also acknowledges with thanks for financial support extended by DST through the J.C. Bose Fellowship. We thank Prof. Franz Lemmermeyer for some valuable suggestions and Paramjit Singh for helping us in the computational part. Finally, our sincere thanks to the anonymous referee for suggesting changes at several places. This work is part of the Ph D thesis of the second author.

\end{document}